\documentclass[12pt]{amsart}
\usepackage{amsfonts,amssymb,amsmath}
\usepackage{srcltx}
\begin{document}

\title[Groups of triangular automorphisms...]{Groups of triangular automorphisms of
a free associative algebra and a polynomial algebra}

\author[Bardakov]{Valeriy G. Bardakov}
\address{Sobolev Institute of Mathematics, Novosibirsk State University, Novosibirsk 630090, Russia}
\email{bardakov@math.nsc.ru}

\author[Neshchadim]{Mikhail V. Neshchadim}
\address{Sobolev Institute of Mathematics, Novosibirsk State University, Novosibirsk 630090, Russia}
\email{neshch@mail.ru}

\author[Sosnovskiy]{Yury V. Sosnovsky}
\address{Novosibirsk State Pedagogical University, Novosibirsk State University, Novosibirsk 630126, Russia}
\email{yury\_sosnovsky@mail.ru}

\begin{abstract}
We study a structure of the group of unitriangular automorphisms of a free associative algebra and a polynomial algebra and prove that this
group is a semi direct product of abelian groups. Using this decomposition we describe a structure of the lower central series and the series
of derived subgroups
for the group of unitriangular automorphisms and prove that every element from the derived subgroup is a commutator.
In addition we prove that the group of unitriangular automorphisms of a free associative algebra of rang more than 2 is not linear and describe
some  two-generated subgroups from these group. Also we give a more simple system of
generators for the group of tame automorphisms  than the system from Umirbaev's paper.
\end{abstract}

\thanks{This research was supported by Federal Target Grant "Scientific and
educational personnel of innovation Russia" for 2009-2013 (government
contract No. 02.740.11.0429)
}
\maketitle

\vspace{0.8cm}

\begin{center}

{\bf Introduction}

\end{center}

\vspace{0.5cm}

In this paper we consider a free associative algebra $A_n = K \langle x_1, x_2,$ $\ldots,$ $x_n \rangle$ and a polynomial algebra
$P_n = K [ x_1, x_2, \ldots, x_n ]$ over a field $K$ of characteristic zero. We assume that all algebras have unity.
By symbol $C_n$ we will denote any from algebras: $A_n$ or $P_n$.
It is evident that $A_1 = P_1$. The group of
$K$-automorphisms of these algebras, i.e. automorphisms which fix elements from $K$ will denote $\mathrm{Aut} \, C_n$.

Despite the great interest to these groups their algebraic  structure isn't well known. In particular,
we do not know a non-trivial system of generators for these groups if
 $n \geq 3$, and a system of defining relations. From the results of Shestakov and Umirbaev
\cite{SU}-\cite{U} follows that
$\mathrm{Aut} \, C_3 \not= \mathrm{TAut} \, C_3$, where
$\mathrm{TAut} \, C_n$ is the group of tame automorphisms of $C_n$.

In the Umirbaen's paper  \cite{U} was written a system of generators for group of tame automorphisms $\mathrm{TAut} \, C_n$ and some system
of relations which is a system of defining relations if $n=3$ (if $n > 3$ then we have an open question: is the written system is
a system of defining relations).

Until now  a description of elements of finite order in group
 $\mathrm{TAut} \, C_n$, $n \geq 3$ is unknown. A corresponding question about the description of the set of involutions was formulated
in \cite[вопрос 14.68]{KT}.

A question about linearity  (i.e. about a faithful representation by finite dimensional matrixes over some field) of $\mathrm{TAut} \, C_n$
was studied in the paper of
Roman'kov, Chirkov, Shevelin
 \cite{R}, where was proved that for  $n \geq 4$ these groups are not linear.
Sosnovskii \cite{S}
proved that for  $n \geq 3$ the group $\mathrm{Aut} \, P_n$ is not linear. This result follows from a description
of the hypercentral series for the  subgroup of unitriangular automorphisms of $\mathrm{Aut} \, P_n$.

In the paper of Gupta, Levchuk and Ushakov \cite{L} were studied generators of the group of tame automorphisms for some nilpotent  algebras.
In addition they proved that any automorphism from  $\mathrm{Aut} \, C_n$ is a product of affine automorphism and automorphism
which maps any generator $x_i,$ $i = 1, 2, \ldots, n$ to an element $x_i + f_i$,
where polynomials  $f_i$ do not contain monomials of degree zero and one.

In the present paper we study a structure of the group of unitriangular automorphisms and prove that this group is a semi direct product
of abelian groups. Using this decomposition we describe a structure of the lower central series and the series of derived subgroups
for the group of unitriangular automorphisms and prove that every element from the derived subgroup is a commutator.
In addition we prove that the group of unitriangular automorphisms of a free associative algebra of rang more than 2 is not linear and describe
some  two-generated subgroups from this group. Also we give a more simple system of
generators for the group of tame automorphisms  than the system from Umirbaev's paper.

We thank the participants of the seminar ``Evariste Galois'' for fruitful discutions. Also we thank  Alexey Belov who points on a connection
between the question on conjugation of any involution from  $\mathrm{Aut} \, P_n$
to diagonal automorphism and Cancelation conjecture.

\vspace{0.8cm}

\begin{center}

{\bf \S~1. Structure of the group of triangular automorphisms}

\end{center}

\vspace{0.5cm}

Let us remember definitions of some automorphisms and subgroups from $\mathrm{Aut} \, C_n$.

For any index  $i\in \left\{1, \ldots, n \right\}$, a constant $\alpha \in K^* = K\backslash \{0\}$ and a polynomial
$f = f(x_1, \ldots , \widehat{x_i}, \ldots ,x_n) \in C_n$ (where the symbol  $\widehat{x}$ denotes that
$x$ does not include in $f$)
{\it the elementary automorphism } $\sigma (i, \alpha, f)$
is an automorphism from $\mathrm{Aut} \, C_n$, which acts on the variables  $x_1, \ldots ,x_n$
by the rule:
$$
\sigma (i, \alpha, f) :
 \left\{
\begin{array}{lcl}
x_i \longmapsto \alpha \, x_i + f, \,\, &    \\
x_j \longmapsto x_j, \,\,           & \mbox{если} & \,\, j \neq i. \\
\end{array}
\right.
$$
The group of the tame automorphisms $\mathrm{TAut} \, C_n$ is generated by all elementary automorphisms.

The group of affine automorphisms $\mathrm{Aff} \, C_n$ is a subgroup of $\mathrm{TAut} \, C_n$ that consists
of automorphisms
$$
x_i\longmapsto a_{i1} x_1+ \ldots + a_{in} x_n + b_i,\,\, i=1, \ldots , n,
$$
where $a_{ij}$, $b_i\in K$, $i,j=1, \ldots ,n$,
and the matrix $(a_{ij})$ is nondegenerated.  The group of affine automorphisms is the semi direct product
$K^n \leftthreetimes \mathrm{GL}_n (K)$  and, in particular, embeds in the group of matrices
 $\mathrm{GL}_{n+1} (K)$.

The group of triangular automorphisms  $T_n = T(C_n)$
of algebra  $C_n$  is generated by automorphisms
$$
x_i\longmapsto \alpha_i x_i+f_i(x_{i+1}, \ldots , x_n),\,\, i=1,\ldots , n,
$$
where $\alpha_i \in K^*$, $f_i\in C_n$ and
$f_n\in K$. If all $\alpha_i =1$ then this automorphism is called
{\it the unitriangular automorphism}. The group of unitriangular automorphism is denoted by
 $U_n = U(C_n)$.
In addition, in $T_n$ singles out the subgroup of diagonal automorphisms $D_n$, which contains the following automorphisms
$$
x_i\longmapsto \alpha_i x_i,\,\, \alpha_i \in K^*,~~~i=1, \ldots , n.
$$

It is easy to see that the following lemma holds

{\bf Lemma 1.} {\it The group $U_n$ is normal in $T_n$ and  $T_n = U_n \leftthreetimes D_n$.}

In the group $U_n$ let us define subgroups
$$
G_i = \mbox{гр}\left( \sigma (i, 1, f) \, | \,  f=f(x_{i+1}, \ldots , x_n)\in C_n   \right), \,\,\, i=1, \ldots , n.
$$
Note that the subgroup $G_i$ is abelian and isomorphic to an additive group of algebra  $C_n$ that is generated by $x_{i+1}, \ldots , x_n$,
$i = 1, \ldots, n-1$ and the subgroup  $G_n $ is isomorphic to the additive group of the field $K$.

The following theorem describes a structure of group of unitriangular automorphisms

{\bf Theorem 1.} { \it  The group $U_n$ is a semi direct product of abelian groups}:
$$
U_n=(\ldots(G_1 \leftthreetimes G_2)\leftthreetimes  \ldots ) \leftthreetimes G_n.
$$

\bigskip

To prove this and the next theorems we will use

{\bf Lemma 2.}
{\it For any non-zero element $a$ of the field  $K$
and any polynomial $g(x_1, \ldots,$ $ x_i,$ $\ldots,$ $x_n) $ from $C_n,$ $n \geq 1,$
there exists a polynomial $f(x_1, \ldots, x_i, \ldots, x_n)$ from $C_n$ such that
the following equality holds}
$$
g(x_1, \ldots, x_i, \ldots, x_n) = f(x_1, \ldots, x_i + a, \ldots, x_n) - f(x_1, \ldots, x_i, \ldots, x_n).
$$

{\bf Proof.} It is evident that  sufficient to consider only the case when $g(x_1, \ldots, x_i, \ldots, x_n)$ is a monomial and $i = 1$.
Let us define a partial order $<$ on the set of monomials from $C_n$. Every monomial has the form
$$
m = b_1 \, x_1^{k_1} \, b_2 \, x_1^{k_2} \, \ldots \, b_s \, x_1^{k_s} \,  b_{s+1},   \eqno{(1)}
$$
where $k_i$ are positive integers; monomials $b_j$ do not contain the variable $x_1$. We will call the number $s$ by
{\it syllable length} of $m$ and will denote the vector
 $(k_1, k_2, \ldots, k_s)$ by $\mu(m)$. If we have other monomial:
$$
m' = c_1 \, x_1^{l_1} \, c_2 \, x_1^{l_2} \, \ldots \, c_r \, x_1^{l_r} \,  c_{r+1}
$$
from $C_n$ then we  say that $m$ is {\it lower} than $m'$ and write $m < m'$ if one from the following conditions holds:

1) $\mathrm{deg}_{x_1} m < \mathrm{deg}_{x_1} m'$, where $\mathrm{deg}_{x_1} m = \sum_{i=1}^s k_i$;

2) $\mathrm{deg}_{x_1} m = \mathrm{deg}_{x_1} m'$, but $s < r$, i.e. the syllable length of $m$ is less than the syllable length of $m'$;

3) $\mathrm{deg}_{x_1} m = \mathrm{deg}_{x_1} m'$, $s = r$, but $(k_1, k_2, \ldots, k_s)$ is lower than $(l_1, l_2, \ldots, l_s)$
related to the lexicographical order.

Let $g(x_1, x_2, \ldots, x_n) =  b_1 \, x_1^{k_1} \, b_2 \, x_1^{k_2} \, \ldots \, b_s \, x_1^{k_s} \,  b_{s+1}$. By induction hypothesis
we will assume  that Lemma is true for all polynomials that are a sum of monomials which are lower $g$. Let
$$
f_1(x_1, x_2, \ldots, x_n) = \frac{1}{(k+1) a} \, b_1 \,  x_1^{k_1+1} \, b_2 \, x_1^{k_2} \, \ldots \, b_s \, x_1^{k_s} \,  b_{s+1}.
$$
It is easy to see that
$$
f_1(x_1 + a, x_2, \ldots, x_n) - f_1(x_1, x_2, \ldots, x_n) = g + g_1(x_1, x_2, \ldots, x_n),
$$
where $g_1$ is a sum of monomials which are lower $g$. By По induction hypothesis there exists a polynomial $f_2$ such that
$$
g_1(x_1, x_2, \ldots, x_n) = f_2(x_1 + a, x_2, \ldots, x_n) - f_2(x_1, x_2, \ldots, x_n).
$$
Then
$$
g = f_1(x_1 + a, x_2, \ldots, x_n) - f_2(x_1 + a, x_2, \ldots, x_n) - f_1(x_1, x_2, \ldots, x_n) + f_2(x_1, x_2, \ldots, x_n).
$$
If we define $f = f_1 - f_2$ then we get the required decomposition. Lemma is proved.




\bigskip

Let us start to prove Theorem 1.

Consider subgroups $H_i$, $i=1, \ldots, n$, of  $U_n$ which consist from automorphisms which fix variables
 $x_i, \ldots, x_{n}$. In addition  let $H_{n+1}=U_n$.
It is evident $H_1=1$, i.e. is a trivial subgroup.

Notice at first that $H_i$ is a normal subgroup in  $U_n$.
Indeed,  if  $\varphi \in H_i$, $\psi \in U_n$ then for variable
  $x_j$, $j \geq i$, we have
$$
x_j^{\psi \varphi {\psi}^{-1}} = (x_j + f_j(x_{j+1}, \ldots, x_n))^{\varphi {\psi}^{-1}} =
 (x_j + f_j(x_{j+1}, \ldots, x_n))^{{\psi}^{-1}} = x_j,
$$
where $x_j^{\varphi} = x_j + f_j(x_{j+1}, \ldots, x_n).$

At second check that
$$
 H_{i+1}=H_i \leftthreetimes G_i,\,\,i=1, \ldots, n.
$$
Since $H_i \triangleleft H_{i+1}$ and $H_i\cap G_i=1$, then sufficient to check that
 $H_{i+1}=H_i \cdot G_i.$ Let $\varphi \in H_{i+1}$ and
$$
x_i^{\varphi}= x_i + g(x_{i+1}, \ldots, x_n),~~ x_j^{\varphi}= x_j, j > i.
$$
Then we can write
$$
 \varphi=\overline{\varphi} \cdot \psi,
$$
where $\overline{\varphi} \in G_i$ and
$x_i^{\overline{\varphi}}=x_i+f(x_{i+1}, \ldots, x_n)$,
but $\psi={\overline{\varphi}}^{-1}\varphi $. Since $x_i^{\psi}=x_i$ then $\psi\in H_i$.
Hence $H_{i+1} = H_i \, G_i$, and then
$$
 U_n=H_{n+1}=H_n \leftthreetimes G_n=
 (H_{n-1} \leftthreetimes G_{n-1})\leftthreetimes  G_n= \ldots =
$$
$$
=(\ldots(G_1 \leftthreetimes G_2)\leftthreetimes  \ldots ) \leftthreetimes G_n.
$$
Theorem 1 is proved.

Let's recall  definitions of the lower central series and  the series of commutators of a group. {\it The lower central series} of group $G$ is called
a series
$$
G = \gamma_1 G \geq \gamma_2 G \geq \ldots,
$$
where $\gamma_{i+1} G = [\gamma_i G, G],$ $i = 1, 2, \ldots$.
{\it The series of commutators} of group $G$ is called a series
$$
G = G^{(0)} \geq G^{(1)} \geq G^{(2)} \geq \ldots,
$$
where $G^{(i+1)} = [G^{(i)}, G^{(i)}],$ $i = 0, 1, \ldots$.

The following theorem gives a description of the lower central series and the  series of commutators of $U_n$.

{\bf Theorem 2.} 1) {\it For $n \geq 2$ in $U_n$  the following equalities hold
$$
\gamma_2 U_n = \gamma_3 U_n = \ldots .
$$
In particular, for $n \geq 2$ the group $U_n$ is not nilpotent.}

2) {\it The group $U_n$ is solvable of degree  $n$ and corresponding commutators subgroups have the form}:
$$
 \begin{array}{l}
  U_n= (\ldots(G_1 \leftthreetimes G_2)\leftthreetimes  \ldots ) \leftthreetimes G_{n},\\
  U_n^{'}= (\ldots(G_1 \leftthreetimes G_2)\leftthreetimes  \ldots ) \leftthreetimes G_{n-1},\\
  .........................................\\
  U_n^{(n-1)}= G_1,\\
  U_n^{(n)}= 1.
 \end{array}
$$

{\bf Proof.}
Let us check that for $i < j$ holds
$$
 [G_i, G_j] = G_i.
$$

The commutator of the element $\varphi = \sigma(i, 1, f)$ from $G_i$ and the element $\psi=\sigma(j, 1, h)$ from $G_j$ has the form
$$
 [\varphi , \psi ] :
 \left\{
\begin{array}{lcl}
x_i \longmapsto x_i + f(x_{i+1},\ldots, x_j+h, \ldots, x_n)-f(x_{i+1}, \ldots, x_j, \ldots, x_n), \,\, &    \\
x_k \longmapsto x_k, \,\,           & \mbox{if} & \,\, k\neq i. \\
\end{array}
\right.
$$
Hence for arbitrary  $\sigma(i, 1, g)$ from $G_i$ take any $h$ from $K^*$ and find by Lemma 2 such polynomial $f = f(x_{i+1},\ldots, x_n)$ for
which  the following equality is true
$$
g = f(x_{i+1},\ldots, x_j + h, \ldots, x_n) - f(x_{i+1},\ldots, x_j,\ldots, x_n).
$$
Then $\sigma(i, 1, g) = [\varphi, \psi]$ and hence
$G_i \leq [G_i,G_j]$. The inverse inclusion is evident.

Let us find the lower central series for $U_2$. By Theorem 1 and by proved equality we have
$$
U_2' = \gamma_2 U_2 = [G_1 \leftthreetimes G_2, G_1 \leftthreetimes G_2] = [G_1, G_2] = G_1.
$$
Further
$$
\gamma_3 U_2 = [\gamma_2 U_2, G_2] = [G_1, G_2] = G_1,
$$
and analogously
$$
\gamma_m U_2 = G_1
$$
for all $m > 3.$
For $n > 2$ the proof is similar.

The second part of Theorem follows from Theorem 1 and the description of the commutator subgroups $[G_i, G_j]$.
Theorem 2 is proved.

\bigskip

Note that the solubility of $U_n$  was proved in \cite{R}.

\bigskip

The next theorem gives an answer to the question about the commutator width of $U_n$.

{\bf Theorem 3.} {\it Every element from the commutator subgroup $U_n^{'}$
is a commutator.}

{\bf Proof.}
From Theorem 2 follows that
$$
 U_n^{'}=
 (\ldots(G_1 \leftthreetimes G_2)\leftthreetimes \ldots) \leftthreetimes G_{n-1}.
$$
Let $\omega\in U_n^{'}$ and
$$
x_i^{ \omega} = x_i+g_i(x_{i+1},...,x_n),\,\, i=1,\ldots, n-1,\,\, x_n^{ \omega} = x_n.
$$
We will find  automorphisms ${\varphi}_i$ from $G_i$, $i=1,
\ldots, n,$  such that
$$
 \omega = [{\varphi}_n, {\varphi}_1\cdot \ldots \cdot {\varphi}_{n-1}],
$$
and
$$
 \varphi_i :
 \left\{
\begin{array}{lcl}
x_i \longmapsto x_i + f_i(x_{i+1},...,x_n), \,\, &   \\
x_j \longmapsto x_j, \,\,           & \mbox{if} & \,\, j\neq i. \\
\end{array}
\right.
$$
Let $f_n=1$.
Using the straight calculus we see that the automorphism $\varphi = [{\varphi}_n, {\varphi}_1\cdot \ldots\cdot {\varphi}_{n-1}]$
acts on the variables by the rule
$$
 x_k^{\varphi} =  x_k + f_k(x_{k+1}+f_{k+1},\ldots,x_{n-1}+f_{n-1},x_n)
                -f_k(x_{k+1}+f_{k+1},\ldots,x_{n-1}+f_{n-1},x_n+1),
$$
$k=1,\ldots,n$. Hence we have to solve the system of equations
$$
 g_k(x_{k+1},\ldots,x_n)=
$$
$$
 =f_k(x_{k+1}+f_{k+1},\ldots,x_{n-1}+f_{n-1},x_n)
                -f_k(x_{k+1}+f_{k+1},\ldots,x_{n-1}+f_{n-1},x_n+1),
$$
relative to polynomials $f_k$, $k=1,\ldots,n-1$.

The equation
$$
 g_{n-1}(x_n)=f_{n-1}(x_n)-f_{n-1}(x_n+1)
$$
has a solution $f_{n-1}(x_n)$ by Lemma 2. Let us assume that the polynomials $f_{n-1},\ldots,f_{k+1}$ are found. Making the change of variables:
$$
  \left\{
 \begin{array}{lcl}
 x_{k+1} + f_{k+1}(x_{k+2},\ldots,x_n)=y_{k+1},   \\
 ..................................... \\
 x_{n-1} + f_{n-1}(x_n)=y_{n-1},   \\
 x_n=y_n,\\
 \end{array}
 \right.
$$
in the $k$-th equation of the system we get the equation
$$
 \widetilde{g}_k(y_{k+1},\ldots,y_n) = f_k(y_{k+1},\ldots,y_{n-1},y_n) - f_k(y_{k+1},\ldots,y_{n-1},y_n+1).
$$
By Lemma 2 there exists a polynomial $f_k$ which satisfy this equation.
Hence there exist the automorphisms ${\varphi}_i\in G_i$, $i=1,\ldots,n,$ such that
$$
 \omega = [{\varphi}_n, {\varphi}_1\cdot \ldots\cdot {\varphi}_{n-1}].
$$
Theorem is proved.

{\bf Question 1.} Is there a finitely generated subgroup from $U_n$ that has the commutator width more than one?

{\bf Question 2.} Is it true that the width of $\mathrm{TAut} \, C_n$, $n \geq 3,$ relatively the elementary automorphisms
$\sigma(i, \alpha, f)$ is finite?

We have defined the group of triangular automorphisms
$$
 T_n=\left\{
            x_i\longmapsto \alpha_i x_i+f_i(x_{i+1},\ldots,x_n) \, \vert \, i=1,\ldots,n,\,\, \alpha_i \in K
     \right\}.
$$
We will call this group by group of {\it upper triangular automorphisms}.
Analogously define a group of {\it lower triangular automorphisms}
$$
 \widetilde{T}_n=\left\{
            x_i\longmapsto \alpha_i x_i+g_i(x_1,\ldots,x_{i-1}) \, \vert \, i=1,\ldots,n,\,\, \alpha_i \in K
     \right\}.
$$

{\bf Question 3.} Is it true that $\mathrm{TAut} \, C_n$ is generated by group of upper triangular automorphisms $T_n$ and
group of lower triangular automorphisms $\widetilde{T}_n$? If it is true, is there  a natural number $m = m(n)$ such that
$$
 \mathrm{TAut} \, (C_n)=\underbrace{T_n\widetilde{T}_n \ldots T_n\widetilde{T}_n}_{m},~~~n > 2?
$$

\vspace{0.8cm}

\begin{center}

{\bf \S~2. Examples of  non-linear subgroups of the group of triangular automorphisms}

\end{center}

\vspace{0.5cm}

In this section we will prove that there exist a countable non-linear subgroups of $U_3$.

At first we prove the following statement.

{\bf Proposition 1.} {\it The group $U_n$, $n \geq 2$, does not contain proper subgroups of finite indexes
and in particular is not residually finite.}

{\bf Proof.} If  $U_n$ contains a subgroup of finite index then it contains a normal subgroup of finite index.
Suppose that $H \lhd U_n$ and $|U_n / H| = m > 1$. Then for every generator $\sigma (i, 1, f)$ of
$U_n$ holds  $(\sigma (i, 1, f))^m \in H$, but it is easy to see that elements $(\sigma (i, 1, f))^m = \sigma (i, 1, m f)$
generate the whole group $U_n$.

The next question is connected with this Proposition.

{\bf Question 4.} Is it true that every finitely generated subgroup of $U_n$:

a) is residually finite?

b) is linear?

Using Proposition 1 it is easy to prove that
 $U_n$, $n \geq 3$ is not linear. For this we will use the next result of Mal'cev \cite{M}:
every soluble linear  group contains a subgroup of finite index the commutator subgroup of which is nilpotent.
From Proposition 1 follows that every subgroup of finite index from $U_3$ is equal to $U_3,$ and from Theorem 2 follows that its
commutator subgroup is not nilpotent.

In fact a stronger assertion holds. To prove this assertion we will use some ideas from \cite{R}.

{\bf Theorem 4.} {\it In
 $U_3$ there exists a countable subgroup which is not linear.}

{\bf Proof.}
Consider the following automorphisms from $U_3$:
$$
\varphi_p :
\left\{%
\begin{array}{l}
x_1\longmapsto x_1+x_2^p, \\
x_2\longmapsto x_2, \\
x_3\longmapsto x_3,
\end{array}%
\right.
~~~
\chi :
\left\{%
\begin{array}{l}
x_1\longmapsto x_1, \\
x_2\longmapsto x_2+x_3, \\
x_3\longmapsto x_3,
\end{array}%
\right.
~~~
\psi :
\left\{%
\begin{array}{l}
x_1\longmapsto x_1, \\
x_2\longmapsto x_2, \\
x_3\longmapsto x_3+1,
\end{array}%
\right.
$$
where $p = 1, 2, \ldots$.

Let us show that a subgroup of $G$, that is generated by elements $\varphi_p,$ $p = 1, 2, \ldots$, $\chi, \psi$ is not linear.
By Mal'cev's theorem which was formulated above, every soluble linear group is a finite extension of a normal subgroup commutator subgroup
of which is nilpotent. If  $G$ is linear then for some natural  $l$ and $r$  the following identity is true in $G$:
$$
[\varphi_p^l, \, _r [\chi^l, \psi^l]] = 1,~~~p = 1, 2, \ldots,
$$
where $[a, \, _1 b] = [a,b] = a^{-1}b^{-1}ab$, $[a, \, _r b] = [[a, \, _{r-1} b], \, b]$.

Using direct calculations we have
$$
[\chi^l, \psi^l] : x_2\longmapsto x_2+l^2,
$$
$$
[\varphi_p^l, [\chi^l, \psi^l]] : x_1\longmapsto x_1+l(x_2+l^2)^p-lx_2^p,
$$
where on unrecorded generators this automorphism acts trivially (we will use this rule of action further).
By induction for arbitrary natural $m$ we get the formula
$$
[\varphi_p^l, \, _m [\chi^l, \psi^l]] : x_1 \longmapsto x_1+\sum\limits_{k=0}^{m} C_m^k \, (-1)^k \, l (x_2+l^2(m-k))^p,
$$
where $C_m^k=\displaystyle\frac{m!}{k!(m-k)!}$.

We must show that for every natural numbers $m$ an $l$ there is a natural $p$, such that the polynomial
$$
\sum\limits_{k=0}^{m} C_m^k \, (-1)^k \, l (x_2+l^2(m-k))^p
$$
is not equal to zero. If we  suppose in this polynomial $x_2=0$ then we get the number
$$
l^{2p+1}\sum\limits_{k=0}^{m}C_m^k \, (-1)^k \, (m-k)^p
$$
and transforming the expression we have
$$
\sum\limits_{k=0}^{m}C_m^k \, (-1)^k \, (m-k)^p=
\sum\limits_{k=0}^{m}C_m^k \, (-1)^k \, \sum\limits_{r=0}^{p} C_p^r \, m^{p-r} \, (-1)^r \, k^r=
$$
$$
 =\sum\limits_{r=0}^{p}
 \left(
  C_p^r \, m^{p-r} \, (-1)^r \sum\limits_{k=0}^{m} C_m^k \, (-1)^k \, k^r
 \right).
$$
To find the sum
$$
  \sum\limits_{k=0}^{m} C_m^k \, (-1)^k \, k^r
$$
consider the function
$$
  f(\xi)=(1-\xi)^m = \sum\limits_{k=0}^{m} C_m^k \, (-1)^k \, \xi^k
$$
with variable $\xi$.
It is obvious that for arbitrary natural $r$ holds
$$
\left(\xi\frac{d}{d\xi}\right)^{r} f(\xi) = \sum\limits_{k=0}^{m} C_m^k \, (-1)^k \, k^r \, \xi^k.
$$
Replace $\xi$ on $e^\zeta$, where $e$ is the base of natural logarithm and let  $\zeta = 0$ we have
$$
  \left.\left(\frac{d}{d\zeta}\right)^{r} (1-e^\zeta)^m\right|_{\zeta=0} =  \sum\limits_{k=0}^{m} C_m^k \, (-1)^k \, k^r.
$$
Hence
$$
 \sum\limits_{r=0}^{p}
 \left(
  C_p^r \, m^{p-r} \, (-1)^r \sum\limits_{k=0}^{m} C_m^k \, (-1)^k \, k^r
 \right)=
$$
$$
 \left.\sum\limits_{r=0}^{p}
 \left(
  C_p^r \, m^{p-r} \, (-1)^r \left(\frac{d}{d\zeta}\right)^{r}
 \right)(1-e^\zeta)^m\right|_{\zeta=0}=
\left.\left(m-\frac{d}{d\zeta}\right)^{p}(1-e^\zeta)^m\right|_{\zeta=0}.
$$
If we take $p=m$ then
$$
 \left.\left(m-\frac{d}{d\zeta}\right)^{p}(1-e^\zeta)^m\right|_{\zeta=0}=
  \left.\sum\limits_{k=0}^{m} C_m^k \, m^{(m-k)} \left(-\frac{d}{d\zeta}\right)^{k} (1-e^\zeta)^m\right|_{\zeta=0}=
$$
$$
  =\left.\left(-\frac{d}{d\zeta}\right)^{m} (1-e^\zeta)^m\right|_{\zeta=0}=m!
$$
Hence for $p = m$ the automorphism $[\varphi_m^l, \, _m [\chi^l, \psi^l]]$ acts on $x_1$ in nontrivial manner.
Since the natural number $p$ can be arbitrary then the theorem is proved.

As was note early the group of affine automorphisms is linear. From this in particular follows that $\mathrm{Aut} \, A_1 = \mathrm{Aut} \, P_1$
is linear.
Thus remains opened the question about linearity of the group $\mathrm{Aut} \, P_2$ which as well known is equal to $\mathrm{Aut} \, A_2$.

\vspace{0.8cm}

\begin{center}

{\bf \S~3. Generators and relations}

\end{center}

\vspace{0.5cm}

Umirbaev \cite{U, U1, U2} defined some system of generators and relations for $\mathrm{TAut} \, C_n$ and for group of automorphisms of
algebra which is free in Shreier-Nilsen variety. In the present section we give some over system of generators and relations.
At first, remember some facts from  \cite{U, U1, U2}. The group $\mathrm{TAut} \, C_n$ is generated by elementary automorphisms
$$
\sigma(i, \alpha, f),~~~\alpha \in K^*, ~~f = f(x_1, \ldots, \widehat{x_i}, \ldots, x_n),~~i = 1, \ldots, n.
$$
We will denote this set by $\mathcal{A}.$ The generators from $\mathcal{A}$ satisfy the next relations which with  correction
on notation of automorphism  $x^{\varphi}$ instead of $\varphi(x)$
have the form
$$
\sigma(i, \alpha, f) \, \sigma(i, \beta, g) = \sigma(i, \alpha \beta, f + \alpha g),  \eqno{(1)}
$$
$$
\tau_{ks} \, \sigma(i, \alpha, f) \, \tau_{ks} = \sigma(i, \alpha, f^{\tau_{ks}}),~~~k \not= i,~~s \not= i, \eqno{(2.1)}
$$
$$
\tau_{is} \, \sigma(i, \alpha, f) \, \tau_{is} = \sigma(s, \alpha, f^{\tau_{is}}), \eqno{(2.2)}
$$
$$
\sigma(i, \alpha, f)^{-1} \sigma(j, \beta, g) \sigma(i, \alpha, f) = \sigma(j, \beta, g^{\sigma(i, \alpha, f)}),
~~~f = f(x_1, \ldots, \widehat{x_i}, \ldots, \widehat{x_j}, \ldots, x_n), \eqno{(3)}
$$
where $\tau_{ks} = \sigma(s, 1, x_k) \sigma(k, 1, -x_s) \sigma(s, -1, x_k)$ is an automorphism of transposition which permutes
the variables $x_k$ and $x_s$. For $n=3$ Umirbaev proved  that the system of relations
(1)--(3) is a system of defining relations. Using this assertion he proved that the kernel of the homomorphism
$\mathrm{TAut} \, A_3 \longrightarrow \mathrm{TAut} \, P_3$ is generated as a normal subgroup by automorphisms
 $\sigma(i, 1, f),$ where $i = 1, \ldots, n,$ and $f$ is an arbitrary element from the commutator ideal of $A_3.$

Define an other set of generators of $\mathrm{TAut} \, C_n$ which similar to the set of generators introdused  by Cohn \cite{Co}
for the group of automorphisms of free Lie algebra, and also Czerniakiewicz \cite{C} for $\mathrm{Aut} \, A_2$. Let
$\varphi(\alpha, f) = \sigma(1, \alpha, f)$, and by symbol
$\tau_{ks},$ $1 \leq k \not= s \leq n,$ we will denote an automorphism that permutes
generator $x_k$ and $x_s$ and fixes other generators:
$$
\tau_{ks} :
 \left\{
\begin{array}{lcl}
x_k \longmapsto x_s, \,\, &   \\
x_s \longmapsto x_k, \,\, &   \\
x_j \longmapsto x_j, \,\,           & \mbox{if} & \,\, j\neq k, s. \\
\end{array}
\right.
$$

In these notations holds

{\bf Theorem 5.} {\it The group $\mathrm{TAut} \, C_n$, $n \geq 3$ is generated by the automorphisms
$$
\varphi(\alpha, f),~~\tau_{ks}, ~~~\alpha \in K^*, ~~f = f(x_2, \ldots, x_n),~~1 \leq k \not= s \leq n,
$$
which satisfy the following relations:
$$
\tau_{ks}^2 = 1,~~\tau_{ks} \, \tau_{lm} = \tau_{lm} \, \tau_{ks}, ~~\tau_{ks} \, \tau_{lk} \, \tau_{ks} = \tau_{ls}, ~~\{ k, s \} \cap \{ l, m \}
= \emptyset, \eqno{(4)}
$$
$$
\varphi(\alpha, f) \varphi(\beta, g) = \varphi(\alpha \beta, \alpha g + f), \eqno{(5)}
$$
$$
\tau_{ks} \, \varphi(\alpha, g) \tau_{ks} = \varphi(\alpha, g^{\tau_{ks}}), ~~k \not=1, s \not= 1, \eqno{(6)}
$$
$$
\tau_{1i} \, \varphi(\alpha, g)^{-1} \, \tau_{1i} \,  \cdot \varphi(\beta, f) \cdot \, \tau_{1i} \, \varphi(\alpha, g) \, \tau_{1i} =
\varphi(\beta, f^{\tau_{1i} \, \varphi(\alpha, g) \tau_{1i}}), ~~g = (x_2, \ldots, \widehat{x_i}, \ldots, x_n). \eqno{(7)}
$$

For $n = 3$ this system of relations is a system of defining relations for} $\mathrm{TAut} \, C_n$.

{\bf Proof.} Note that the generators $\tau_{ks}$ and relations (4) define the symmetric group $S_n$. Denote by $\mathcal{B}$ the system
of generators from Theorem 4. Show that every generator from the set $\mathcal{A}$
express in terms of generators from $\mathcal{B}$. Consider $\tau_{1i} \, \sigma(i, \alpha, f) \, \tau_{1i}.$ Acting by
$\tau_{1i} \, \sigma(i, \alpha, f) \, \tau_{1i}$ on the generators we get
$$
\left.%
\begin{array}{ll}
x_1^{\tau_{1i} \, \sigma(i, \alpha, f) \, \tau_{1i}} = x_i^{\sigma(i, \alpha, f) \, \tau_{1i}} =  (\alpha x_i + f)^{\tau_{1i}} = \alpha x_1 + f^{\tau_{1i}}, &  \\
& \\
x_i^{\tau_{1i} \, \sigma(i, \alpha, f) \, \tau_{1i}} = x_1^{\sigma(i, \alpha, f) \, \tau_{1i}} =  x_1^{\tau_{1i}} = x_i, & \\
& \\
x_j^{\tau_{1i} \, \sigma(i, \alpha, f) \, \tau_{1i}} = x_j^{\sigma(i, \alpha, f) \, \tau_{1i}} =  x_j^{\tau_{1i}} = x_j, ~~~ \mbox{if}~~j \not= 1, i. & \\
& \\
\end{array}%
\right.
$$
Hence
$$
\tau_{1i} \, \sigma(i, \alpha, f) \, \tau_{1i} = \varphi(\alpha, f^{\tau_{1i}}).
$$
From this
$$
\sigma(i, \alpha, f) = \tau_{1i} \, \varphi(\alpha, f^{\tau_{1i}}) \, \tau_{1i}. \eqno{(8)}
$$
Therefore any element from $\mathcal{A}$ expresses in terms of generators from  $\mathcal{B}$.

Rewrite the relations (1)--(3) in generators from $\mathcal{B}$. Relation (1) takes the following form
$$
\tau_{1i} \, \varphi(\alpha, f^{\tau_{1i}}) \, \tau_{1i} \cdot \tau_{1i} \, \varphi(\beta, g^{\tau_{1i}}) \, \tau_{1i} =
\tau_{1i} \, \varphi(\alpha \beta, f^{\tau_{1i}} + \alpha g^{\tau_{1i}}) \, \tau_{1i}.
$$
Conjugating both  sides by automorphism $\tau_{1i}$ we get
$$
 \varphi(\alpha, f^{\tau_{1i}})  \cdot \varphi(\beta, g^{\tau_{1i}})  =
\varphi(\alpha \beta,  f^{\tau_{1i}} + \alpha g^{\tau_{1i}}),
$$
but this is the relation (5).

Relation (2.1) for $i = 1$  has the form
$$
\tau_{ks}  \, \sigma(1, \alpha, f)  \, \tau_{ks} =
\sigma(1, \alpha,  f^{\tau_{ks}}),~~~k \not= 1, s \not= 1
$$
or
$$
\tau_{ks} \,  \varphi(\alpha, f) \, \tau_{ks} =
\varphi(\alpha,  f^{\tau_{ks}}),
$$
but this is the relation (6).

Relation (2.2) in new generators has the form
$$
\tau_{is} \, \tau_{1i} \, \varphi(\alpha, f^{\tau_{1i}}) \, \tau_{1i} \, \tau_{is} =
\tau_{1s} \, \varphi(\alpha, f^{\tau_{is} \, \tau_{1s}}) \, \tau_{1s}.
$$
Conjugating both sides by automorphism $\tau_{1s}$ we get
$$
\tau_{1s} \, \tau_{is} \, \tau_{1i} \, \varphi(\alpha, f^{\tau_{1i}}) \, \tau_{1i} \, \tau_{is} \, \tau_{1s} =
\varphi(\alpha, f^{\tau_{is} \, \tau_{1s}}).
$$
From (4) follows the equality
$$
\tau_{1s} \, \tau_{is} \, \tau_{1i} = \tau_{is},
$$
using which we have
$$
\tau_{is} \,  \varphi(\alpha, f^{\tau_{1i}}) \,  \tau_{is} =
\varphi(\alpha, f^{\tau_{is} \, \tau_{1s}}).
$$
Let
$$
f^{\tau_{1i}} = g(x_1, \ldots, \widehat{x_i}, \ldots, x_n)
$$
and from the equality $\tau_{is} \, \tau_{1s} = \tau_{1i} \, \tau_{is}$ we can write:
$f^{\tau_{is} \, \tau_{1s}} = f^{\tau_{1i} \, \tau_{is}} = g^{\tau_{is}}$ and this relation is reduced to
$$
\tau_{is} \,  \varphi(\alpha, g) \,  \tau_{is} =
\varphi(\alpha, g^{\tau_{is}}),
$$
and this is (6).

Relation (3) in new generators has the form
$$
\tau_{1i} \, \varphi(\alpha, f^{\tau_{1i}})^{-1} \, \tau_{1i} \cdot  \tau_{1j} \, \varphi(\beta, g^{\tau_{1j}}) \, \tau_{1j} \cdot
\tau_{1i} \, \varphi(\alpha, f^{\tau_{1i}}) \, \tau_{1i} =
\tau_{1j} \, \varphi(\beta, g^{\tau_{1i} \, \varphi(\alpha, f^{\tau_{1i}}) \, \tau_{1i} \, \tau_{1j}}) \, \tau_{1j}.
$$
Conjugating both sides by automorphism $\tau_{1j}$ we get
$$
\tau_{1j} \, \tau_{1i} \, \varphi(\alpha, f^{\tau_{1i}})^{-1} \, \tau_{1i} \, \tau_{1j} \cdot   \varphi(\beta, g^{\tau_{1j}}) \,  \cdot
\tau_{1j} \, \tau_{1i} \, \varphi(\alpha, f^{\tau_{1i}}) \, \tau_{1i} \, \tau_{1j} =
\varphi(\beta, g^{\tau_{1i} \, \varphi(\alpha, f^{\tau_{1i}}) \, \tau_{1i} \, \tau_{1j}}). \eqno{(8)}
$$
Since $\tau_{1j} \, \tau_{1i} = \tau_{1i} \, \tau_{ij}$ then
$$
\tau_{1j} \, \tau_{1i} \, \varphi(\alpha, f^{\tau_{1i}}) \, \tau_{1i} \, \tau_{1j} = \tau_{1i} \, \tau_{ij} \, \varphi(\alpha, f^{\tau_{1i}})
 \, \tau_{ij} \, \tau_{1i}
= \tau_{1i} \, \varphi(\alpha, f^{\tau_{1i} \, \tau_{ij}}) \, \tau_{1i}, \eqno{(9)}
$$
where we used the proven relation (6).
Using (9) and note that the element $f^{\tau_{1i}} = h(x_2, \ldots, \widehat{x_j}, \ldots, x_n)$
does not contain $x_1$ and $x_j$ we can rewrite (8) in the form
$$
\tau_{1i} \, \varphi(\alpha, h^{\tau_{ij}})^{-1} \, \tau_{1i} \cdot   \varphi(\beta, g^{\tau_{1j}}) \,  \cdot
\tau_{1i} \,  \varphi(\alpha, h^{\tau_{ij}}) \, \tau_{1i}  =
\varphi(\beta, g^{\tau_{1i} \, \varphi(\alpha, h) \, \tau_{1i} \, \tau_{1j}}). \eqno{(10)}
$$
Let $h^{\tau_{ij}} = l(x_2, \ldots, \widehat{x_i}, \ldots, x_n),$ $g^{\tau_{1j}} = m(x_2, \ldots, x_n)$
and take into account that $l^{\tau_{ij}} = h$, $m^{\tau_{1j}} = g$ rewrite the equality (10) in the form
$$
\tau_{1i} \, \varphi(\alpha, l)^{-1} \, \tau_{1i} \cdot   \varphi(\beta, m) \,  \cdot
\tau_{1i} \,  \varphi(\alpha, l) \, \tau_{1i}  =
\varphi(\beta, m^{\tau_{1j} \, \tau_{1i} \, \varphi(\alpha, l^{\tau_{ij}}) \, \tau_{1i} \, \tau_{1j}}).
$$
Using the equality  $\tau_{1j} \, \tau_{1i} = \tau_{1i} \, \tau_{ij}$ we can write
$$
\tau_{1j} \, \tau_{1i} \, \varphi(\alpha, l^{\tau_{ij}}) \, \tau_{1i} \, \tau_{1j} =
\tau_{1i} \, \tau_{ij} \, \varphi(\alpha, l^{\tau_{ij}}) \, \tau_{ij} \, \tau_{1i} =
\tau_{1i} \, \varphi(\alpha, l) \, \tau_{1i},
$$
and finally we have
$$
\tau_{1i} \, \varphi(\alpha, l)^{-1} \, \tau_{1i} \cdot   \varphi(\beta, m) \,  \cdot
\tau_{1i} \,  \varphi(\alpha, l) \, \tau_{1i}  =
\varphi(\beta, m^{ \tau_{1i} \, \varphi(\alpha, l) \, \tau_{1i}})
$$
but this is the relation (7).

{\bf Question 5.} What automorphisms we have to add to the set of the tame automorphisms in order that this set generates $\mathrm{Aut} \, C_3$?

{\bf Question 6.} Is it possible to construct a normal form for words in  $\mathrm{TAut} \, C_3$?

\vspace{0.8cm}

\begin{center}

{\bf \S~4. Some groups generated by two elements}

\end{center}

\vspace{0.5cm}

In this section we will consider subgroups which are generated by two elementary automorphisms. The full description of all these subgroups
does not know even for two dimension linear groups  (see the survey \cite{Mer}). In particular, we do not known when two
elementary transvections  generat a free subgroup.

Consider 2-generated subgroups: $\langle \sigma (i, \alpha, f), \sigma (j, \beta, g) \rangle$ of $\mathrm{Aut} \, C_n$.
Depends weather  index $i$ equal to $j$ or not, we will consider two cases, and without  loss of generality we can assume
that  $i = j = 1$ or $i = 1,$ $j = 2$.
The following Theorem is valid.

{\bf Theorem 6.} {\it Let $\varphi = \sigma (1, \alpha, f)$, $\psi = \sigma (2, \beta, g)$, $\mathrm{deg}_{x_2} f = p$, $\mathrm{deg}_{x_1} g = q$.
Then for $p \cdot q \geq 2$ the following isomorphism holds}
$$
\langle \varphi, \psi \rangle \simeq \langle \varphi \rangle * \langle \psi \rangle.
$$

{\bf Proof.} Let $w$ be a non-trivial word in the free product $\langle \varphi \rangle * \langle \psi \rangle$. Let us show that $w$
is representing a non-trivial element from $\langle \varphi, \psi \rangle$.
We can assume that
$$
w \not\in \langle \varphi \rangle,~~~w \not\in \langle \psi \rangle
$$
and using a conjugation if needed we can obtain that $w$ has the form
$$
w = \varphi^{k_1} \, \psi^{l_1} \, \ldots \, \varphi^{k_m} \, \psi^{l_m},~~~k_i, l_i \in \mathbb{Z},
$$
where $0 < k_i < |\varphi|$ if the order $|\varphi|$ of $\varphi$ is finite and $k_i \not= 0$ if $|\varphi| = \infty$; analogously
$0 < l_i < |\psi|$ if the order  $|\psi|$ of $\psi$ is finite and $l_i \not= 0$ if $|\psi| = \infty$, $i = 1, \ldots, m$.

Let a polynomial $h$ from $C_n$ have the form
$$
h = x_1^s \, A_0(x) + B_0(x),
$$
where  $s \geq 1$, the polynomial  $A_0(x)$ does not contain the variables $x_1,$ $x_2,$ and total degree of the $B_0(x)$ by the variables
$x_1,$ $x_2$ is strictly less than  $s$. Then
$$
h^{\varphi^{k_1} \psi^{l_1}} =  \left( (\alpha^{k_1} \, x_1 + (1 + \alpha + \ldots + \alpha^{k_1-1}) \cdot
f(\beta^{l_1} \, x_2 + (1 + \beta + \ldots + \beta^{l_1 - 1}) \, g, x_3, \ldots, x_n)\right)^s \cdot 
$$
$$
\cdot A_0(x) + B_0 \left( \alpha^{k_1} \, x_1 + (1 + \alpha + \ldots + \alpha^{k_1-1}) \, f(\beta^{l_1} \, x_2 + (1 + \beta + \ldots + \beta^{l_1 - 1}) \, g, x_3, \ldots, x_n), \right.
$$
$$
\left. \beta^{l_1} \, x_2 + (1 + \beta + \ldots +
\beta^{l_1 - 1}) \, g, x_3, \ldots, x_n\right).
$$
From the conditions on $A_0(x)$, $B_0(x)$ and   form of the left side we have
$$
h^{\varphi^{k_1} \psi^{l_1}} = x_1^{s p q} \, A_1(x) + B_1(x),
$$
where  $A_1(x)$ does not contain $x_1,$ $x_2,$ and  degree of $B_1(x)$ by the variables $x_1,$ $x_2$ is strictly less than $s p q$.
By induction of the syllable length $m$ of the word $w$ we get that
$$
h^{w} = x_1^{s (p q)^m} \, A_m(x) + B_m(x),
$$
where $A_m(x)$ does not contain $x_1,$ $x_2,$ and degree of the $B_m(x)$ by the variables $x_1,$ $x_2$ is strictly less than $s (p q)^m$.
In particular, if $h = x_1$ then $x_1^{w} \not= x_1$. Hence $w$ is a non-trivial element from the group $\langle \varphi, \psi \rangle$.
For the free associative  algebra  $A_n$ we must use the homomorphism
$A_n \longrightarrow P_n$. Theorem is proved.

We do not know the structure of $\langle \varphi, \psi \rangle$ in the case of $p q < 2$.
For example, if $\varphi = \sigma (1, 1, \mu x_2)$, $\psi = \sigma (2, 1, \mu x_1)$,  $\mu \in \mathbb{C},$
then $\langle \varphi, \psi \rangle$ is isomorphic to the group which is generated by two transvections
$$
\left(%
\begin{array}{cc}
  1 & \mu \\
0 & 1 \\
\end{array}%
\right),~~~
\left(%
\begin{array}{cc}
  1 & 0 \\
\mu & 1 \\
\end{array}%
\right),
$$
but for this group we do not know the conditions of $\mu$ under which this group is free (see \cite{Mer}).

We will describe some more cases when we can know the structure of $\langle \varphi, \psi \rangle$. Holds

{\bf Proposition 2.} {\it Let $\varphi = \sigma (1, \alpha, f)$, $\psi = \sigma (1, \beta, g)$ where  $f\not= 0$ and $g\not= 0$. Then}

1) {\it if $\alpha \not= 1$ or $\beta \not= 1$ then the group $\langle \varphi, \psi \rangle$ is metabelian};

2) {\it if $\alpha = \beta = 1$ then}
$$
\langle \varphi, \psi \rangle
\simeq
 \left\{
\begin{array}{lll}
 \mathbb{Z}\times  \mathbb{Z},   \,\, & \mbox{if} &
           \,\, \frac{f}{g} \not\in \mathbb{Q},    \\
           & & \\
 \mathbb{Z}, \,\, & \mbox{if} &
           \,\, \frac{f}{g} \in \mathbb{Q}. \\
\end{array}
\right.
$$

{\bf Proof.} If $\alpha \not= 1$ or $\beta \not= 1$ then the commutator $[\varphi, \psi]$ has the form
$$
[\varphi, \psi] = \sigma (1, 1, \beta^{-1} (1 - \alpha^{-1}) g - \alpha^{-1} (1 - \beta^{-1}) f)
$$
and lies in the abelian group $G_1$. Hence the group  $\langle \varphi, \psi \rangle$ is metabelian.

If $\alpha = \beta = 1$ then $a, b \in G_1$ and
$$
x_1^{\varphi^k \, \psi^l} = x_1 + k \, f + l \, g.
$$
For $f / g = p / q$ where $p, q \in \mathbb{Z},$ the group $\langle \varphi, \psi \rangle$ is isomorphic to the subgroup $\langle 1, p / q \rangle$
of the additive group of rational numbers. It is easy to see that this group is infinite cyclic group.

For  $f / g \not\in \mathbb{Q}$  the inequality $k \, f + l \, g \not= 0$ is true for all integer numbers $k, l$ which are not both equal to zero.
Hence $\langle \varphi, \psi \rangle \simeq  \mathbb{Z}\times  \mathbb{Z}$.
This completes the proof.

{\bf Proposition 3.} {\it Let $\varphi = \sigma (1, \alpha, f)$, $\psi = \sigma (2, \beta, g)$ where  $g \in K$.
Then the group $\langle \varphi, \psi \rangle$ is a metabelian group.}

{\bf Proof.} Since the commutator $[\varphi, \psi]$ has the form
 $$
[\varphi, \psi] = \sigma (1, 1, \alpha^{-1} \, f (\beta \, x_2 + g, x_3, \ldots, x_n) - \alpha^{-1} \, f (x_2, x_3, \ldots, x_n)),
 $$
then the derived subgroup of $\langle \varphi, \psi \rangle$ lies in the abelian group $G_1$. Proposition is proved.

 {\bf Question 7.} Is it possible that a group, which is generated by two elementary automorphisms  $\sigma (1, \alpha, f)$, $\sigma (2, \beta, g)$,
is a solvable  of degree more than 3?

{\it The poisson group}
$$
\mathcal{H} = \langle F_2 \times F_2, t ~||~t (g, g) t^{-1} = (g, 1),~~g \in F_2 \rangle
$$
was defined in \cite{FP}. This group is not linear. Also in \cite{FP} was proved that the group of automorphisms
$\mathrm{Aut} \, F_3$ of a  free group $F_3$ contains the poisson group.

 {\bf Question 8.} Does the group $\mathrm{Aut} \, C_n$ contain the poisson group?

 {\bf Question 9.} Does the  Tits alternative hold in $\mathrm{Aut} \, C_n$ or in  $\mathrm{TAut} \, C_n$?


\vspace{0.8cm}

\begin{center}

{\bf \S~5. Elements of finite order in $\mathrm{Aut} \, C_n$}

\end{center}

\vspace{0.5cm}

Evidently that   $\mathrm{Aut} \, C_n$  contains as subgroup $\mathrm{GL}_n(K)$ for $n \geq 2$ and hence contains elements of finite order.
It is easy to describe elementary automorphisms which have finite orders.

{\bf Proposition 4.} {\it The automorphism $\sigma(i, \alpha, f)$ of $C_n$ has an order $m > 1$ if and only if the element
$\alpha$ has the order $m$ in the multiplicative group of the field $K$.}

{\bf Proof.} If the automorphism $\sigma(i, \alpha, f)$ has the order $m_1 > 1$ then from the equality
$$
x_i = x_i^{\sigma^{m_1}(i, \alpha, f)} = \alpha^{m_1} x_i + (1 + \alpha + \ldots + \alpha^{m_1-1}) f
$$
follows that $\alpha \not= 1,$ $\alpha^{m_1} = 1$, i.~e. $\alpha$ is an element from $K^*$ of finite order $m_2$ and $1 < m_2 \leq m_1$.
If $\alpha$ is an element from $K^*$ of the order $m_2$
and $m_2 > 1$ then $\alpha^{m_2} = 1$, $1 + \alpha + \ldots + \alpha^{m_2-1} = 0$. Hence $\sigma^{m_2}(i, \alpha, f) = 1$, i.~e. the automorphism
$\sigma(i, \alpha, f)$ has the finite order $m_1$ and $1 < m_1 \leq m_2$. Hence $m_1 = m_2$ and Proposition is proved.

We describe some subgroups of $\mathrm{Aut} \, C_n$ that do not have a torsion.
We will denote by $IA(C_n)$ a subgroup of $\mathrm{Aut} \, C_n$ that is generated by automorphisms that are identical by module
of an ideal $R^2$ where $R$ is an ideal that is generated by polynomials without free term. These subgroups are studied
in the paper \cite{L}. More generally, for any natural number $k \geq 1$ we
define a subgroup $IA^{(k)}(C_n)$ that is generated by automorphisms that are identical by module of the ideal $R^{k+1}$. It is evident that
$$
IA(C_n) = IA^{(1)}(C_n) \leq IA^{(2)}(C_n) \leq IA^{(3)}(C_n) \leq \ldots
$$

{\bf Proposition 5.} {\it The group  $IA(C_n)$ does not have a torsion.}

{\bf Proof.} Consider the case of $IA(P_n).$ Let $\varphi$ be an automorphism from $IA(P_n)$. Find the minimal natural number $k$ such that
$\varphi$ is identical by the module $R^k$ but is not identical by the module $R^{k+1}$. Then by the modular $R^{k+1}$ the automorphism
$\varphi$ has the form
$$
\overline{\varphi} : x_i \longmapsto x_i + f_i (x_1, x_2, \ldots, x_n),~~~i = 1, 2, \ldots, n,
$$
where $f_i (x_1, x_2, \ldots, x_n)$ is a homogeneous polynomial of degree $k$ or is equal to zero. Then by the module  $R^{k+1}$ we have
$$
\overline{\varphi}^m : x_i \longmapsto x_i + m f_i (x_1, x_2, \ldots, x_n),~~~i = 1, 2, \ldots, n,
$$
and it is evident that for every non-zero natural number $m$ the automorphism $\overline{\varphi}^m$ is not identical. Proposition is proved.

If we denote by $TIA(C_n)$ the intersection of $IA(C_n)$ and $TA(C_n)$ then we can formulate

{\bf Question 10.} Is it true  that $TIA(C_n)$ is generated by the automorphisms $\sigma(i, 1, f)$,
$f = f(x_1, \ldots, \widehat{x_i}, \ldots, x_n) \in R^2$?

The following question is well known  (see for example \cite{K}).

{\bf Question 11.} Is it true that any element of finite order from $\mathrm{Aut} \, C_n$ conjugates in $\mathrm{Aut} \, C_n$
to some element from $\mathrm{GL}_n(K)$?

For involutions this question is formulated in \cite[question 14.68]{KT}.

For the elementary automorphisms hold

{\bf Theorem 7.} {\it Every non trivial automorphism $\sigma(i, \alpha, f)$ of $C_n$
conjugates to some diagonal automorphism if and only if $\alpha \not= 1$.}

{\bf Proof.} We will prove this theorem for the case $i = 1$. If $\alpha \not= 1$ then from the equality
$$
x_1^{\sigma^{-1}(1, 1, - (\alpha - 1)^{-1} f) \, \sigma(1, \alpha, f) \, \sigma(1, 1, - (\alpha - 1)^{-1} f)} = \alpha x_1
$$
follows that $\sigma(1, \alpha, f)$ conjugates to a diagonal automorphism by elementary automorphism
 $\sigma(1, 1, - (\alpha - 1)^{-1} f)$.

Let $\alpha = 1$. Assume that the following equality holds
$$
\varphi^{-1} \sigma(1, 1, f) \varphi = \delta,
$$
where $\delta$ is a diagonal automorphism and $\varphi$ is some automorphism from $\mathrm{Aut} \, C_n$. Denote by
$$
y_1 = x_1^{\varphi},~~y_2 = x_2^{\varphi},~~\ldots, ~~y_n = x_n^{\varphi}.
$$
Evidently that the polynomials  $y_1, y_2, \ldots, y_n$ generate $C_n$. Rewriting our equality in the form
$$
 \sigma(1, 1, f) \varphi = \varphi \delta,
$$
and acting by the both sides on the generators $x_1, x_2, \ldots, x_n$ we get the equalities
$$
x_1^{\sigma (1, 1, f) \varphi} = (x_1 + f(x_2, \ldots, x_n))^{\varphi} = y_1 + f(y_2, \ldots, y_n) = x_1^{\varphi \delta} = y_1^{\delta},
$$
$$
x_2^{\sigma (1, 1, f) \varphi} = x_2^{\varphi} = y_2 = x_2^{\varphi \delta} = y_2^{\delta},~~\ldots,~~y_n = y_n^{\delta}.
$$
Rewriting our equality in the form
$$
y_1^{\delta} - y_1 = f(y_2, \ldots, y_n),
$$
and acting on the both sides by automorphism $\delta$ we have
$$
y_1^{\delta^2} - y_1^{\delta} = f(y_2^{\delta}, \ldots, y_n^{\delta}) = f(y_2, \ldots, y_n).
$$
Then
$$
(y_1^{\delta^2} - y_1^{\delta}) - (y_1^{\delta} - y_1) = y_1^{\delta^2} - 2 y_1^{\delta} + y_1 = 0.
$$
Since $\delta$ acts on  each monomial of the polynomial $y_1$ by dilatation then from the last equality follows
that all coefficients of dilatation satisfy to the equation
$$
c^2 - 2 c + 1 = 0,
$$
i.e. all coefficients of dilation are equal to 1. Hence, $y_1^{\delta} = y_1$ and $\delta$ is the identical automorphism and then
$\sigma(1, 1, f) = 1$ but it is contradiction. Theorem is proved.

It is readily  seen that this theorem is true for the fields of non zero characteristic.

Let's notice that if involution from $\mathrm{Aut} \, C_n$ conjugates to a linear automorphism then it conjugate to some diagonal automorphism.
Known that from the positive answer to the question about involutions follows the positive solution of Cancellation conjecture
(see for example \cite{MSJ}).
It is quite possible that this fact well known to specialists by polynomial maps but we did not find the proof in the literature.
We will prove this fact following to idea of A.~Ya.~Belov.

{\bf Cancellation conjecture.} Let $R$ be a finitely generated commutative algebra over a field $K$ and
$R[z] \simeq K[x_1, x_2, \ldots, x_{n-1}, y]$. Is it true that then
$R$ is isomorphic to $K[x_1, x_2,$ $\ldots,$ $x_{n-1}]$?

The following statement holds.

{\bf Theorem 8.} {\it From the fact that any involution from $\mathrm{Aut} \, P_n$ conjugates to some linear automorphism follows
the positive solution of Cancellation conjecture.}

{\bf Proof.} Let the following isomorphism holds
$$
R[z] \simeq K[x_1, x_2, \ldots, x_{n-1}, y].
$$
Without loss of generality we will assume that these algebras are not only isomorphic but  equal.
This allows  us to use two different notations for elements. Denote the set of variables $x_1, x_2, \ldots, x_{n-1}$ by $x$.
Consider an involution  $\varphi$ of $R[z]$ that is defined by the rule:
$$
\varphi :
\left\{%
\begin{array}{lc}
  z \longmapsto -z, &  \\
  r \longmapsto r, & r \in R. \\
\end{array}%
\right.
$$
By assumption the involution $\varphi$ of $K[x, y]$ conjugates to some linear automorphism.
Then we can choose polynomials $f_1(x, y)$,
$f_2(x, y)$, $\ldots$, $f_{n-1}(x, y)$, $g(x, y)$ from  $K[x, y]$, on which $\varphi$ acts by the following manner:
$$
\varphi :
\left\{%
\begin{array}{ll}
  f_i \longmapsto \varepsilon_i \, f_i, & \varepsilon_i = \pm 1,~~i = 1, 2, \ldots, n-1,  \\
  g \longmapsto \varepsilon \, g, & \varepsilon = \pm 1,  \\
\end{array}%
\right.
$$
and
$$
K[x_1, x_2, \ldots, x_{n-1}, y] = K[f_1, f_2, \ldots, f_{n-1}, g].
$$
Without loss of generality we can assume that
$$
\varphi :
\left\{%
\begin{array}{ll}
  x_i \longmapsto \varepsilon_i \, x_i, & \varepsilon_i = \pm 1,~~i = 1, 2, \ldots, n-1,  \\
  y \longmapsto \varepsilon \, y, & \varepsilon = \pm 1.  \\
\end{array}%
\right.
$$
Note that among the coefficients $\varepsilon_1, \varepsilon_2, \ldots, \varepsilon_{n-1}, \varepsilon $ only one is equal to $-1$.
Indeed if we consider factor algebra of $K[x, y]$ by ideal $I$ that is generated by elements
$$
x_i \, x_j, ~~x_i \, y,~~y^2,~~~i, j = 1, 2, \ldots, n-1,
$$
then we get an algebra which is isomorphic to the vector space $K^{n+1}$, and the automorphism $\varphi$ induces  a linear operator with eigenvalues $\varepsilon_1,$ $\varepsilon_2,$ $\ldots,$
$\varepsilon_{n-1},$ $\varepsilon $, $1$. On the other hand the automorphism $\varphi$ induces on the factor algebra $R[z] / I$
a linear operator with the eigenvalues
$$
-1, 1, 1, \ldots, 1,
$$
since $z^{\varphi} = -z,$ $r^{\varphi} = r$, $r \in R$.

Hence we can assume that $\varphi$ acts by the following manner:
$$
x_1^{\varphi} = x_1,~~ x_2^{\varphi} = x_2,~~ \ldots,~~ x_{n-1}^{\varphi} = x_{n-1},~~ y^{\varphi} = -y.
$$
Write the element $z \in K[x, y]$ in the form
$$
z = a_0(x) + a_1(x) \, y + \ldots + a_k(x) \, y^k,
$$
where  $a_i(x) \in K[x]$ and act by $\varphi$ we get the equality
$$
z^{\varphi} = -z =  a_0(x) - a_1(x) \, y + \ldots + (-1)^k \, a_k(x) \, y^k.
$$
Hence $z$ contains only odd powers of $y$ and we can write
$$
z = y \, p(x, y^2)
$$
for some polynomial $p(x, y^2)$ which depends on $x, y^2$.

On the other hand we can decompose the element $y$ from $R[z]$ by the powers of $z$:
$$
y = b_0(r) + b_1(r) \, z + \ldots + b_l(r) \, z^l,
$$
where $b_i(r) \in R$. Taking into account the equality
$$
y^{\varphi} = -y,~~z^{\varphi} = -z,~~r^{\varphi} = r,~~r \in R,
$$
as well as above we get
$$
y = z \, q(z^2, r),
$$
for some polynomial $q(z^2, r)$ from $R[z]$. Hence
$$
z = y \, p(x, y^2) = z \, q(z^2, r) \, p(x, y^2),
$$
and from this
$$
1 = q(z^2, r) \, p(x, y^2).
$$
Last equality means that $p(x, y^2)$ is an invertible in  $K[x, y]$ and  $p(x, y^2) = p \in K \setminus \{ 0 \}$. Then
$z = y \, p$ and
$$
R[z] = R[y] = K[x, y].
$$
If we consider a factor algebra by an ideal $J$ that is generated by element $y$ and take into account that $R \cap J = 0$ we get
 $R = K[x]$. Theorem is proved.

Let us formulate a hypothesis from which  follows the diagonalization of involutions. Let $\varphi$ be an involution from
$\mathrm{Aut} \, P_n$, i.e. $\varphi^2 = 1$. For an arbitrary polynomial $f$ from $P_n$ holds
$$
(f + f^{\varphi})^{\varphi} = f + f^{\varphi},~~~(f - f^{\varphi})^{\varphi} = -(f - f^{\varphi}).
$$
Hence we can write  $f$ in the form
$$
f = f_1 + f_2,~~f_1^{\varphi} = f_1,~~f_2^{\varphi} = -f_2,
$$
where
$$
f_1 = \frac{f + f^{\varphi}}{2},~~~f_2 = \frac{f - f^{\varphi}}{2}.
$$
Therefore the space  $P_n$ is the direct sum
$$
P_n = \mathrm{Fix} \, \varphi \oplus \mathrm{IFix} \, \varphi,
$$
where $\mathrm{Fix} \, \varphi$ is the set of polynomials that are fixed by the action of $\varphi$, and
$\mathrm{IFix} \, \varphi$ is the set of polynomials that reverse sign to the opposite
by the action of $\varphi$.
Write the generator $x_i$ in the form
$$
x_i = y_i + z_i,~~y_i^{\varphi} = y_i,~~z_i^{\varphi} = -z_i,~~i = 1, 2, \ldots, n,
$$
and assume that holds

{\bf Hypothesis.} From the set of polynomials $y_i, z_i$, $i = 1, 2, \ldots, n$ it is possible to choose a basis of the algebra $P_n$.

Then the automorphism $\varphi$ has in this basis a diagonal matrix. Indeed, let us assume that
$$
K[x_1, \ldots, x_n] = K[y_1, \ldots, y_l, z_{l+1}, \ldots, z_n].
$$
Consider the automorphism $\psi$ that is defined by the following manner:
$$
\psi :
\left\{%
\begin{array}{cc}
  x_i \longmapsto y_i, & 1 \leq i \leq l,\\
  x_i \longmapsto z_i, & l+1 \leq i \leq n.\\
\end{array}%
\right.
$$
It is easy to check that
$$
x_i^{\psi \varphi \psi^{-1}} = x_i,~~~1 \leq i \leq l,
$$
$$
x_j^{\psi \varphi \psi^{-1}} = -x_j,~~~l+1 \leq j \leq n.
$$
Hence $\varphi$ conjugates to a diagonal automorphism.

In the present article we considered the field $K$ of zero characteristic.

{\bf Question 12.} What from the results of the present article are true for fields of non-zero characteristic?

\end{document}